\documentclass[12pt]{article}

\usepackage{a4}
\usepackage{parskip}

\usepackage{amsmath}
\usepackage{amsfonts}
\usepackage{amssymb}
\usepackage{theorem}
\usepackage{latexsym}

\newcommand{\C}{\mathbb{C}}

\newcommand{\N}{\mathbb{N}}\renewcommand{\P}{\mathbb{P}}
\newcommand{\R}{\mathbb{R}}

\newcommand{\cL}{{\cal L}}

\newtheorem{theorem}{Theorem}[section]

\newtheorem{proposition}[theorem]{Proposition}
\newtheorem{lemma}[theorem]{Lemma}

\newtheorem{remark}[theorem]{Remark}
\newtheorem{example}[theorem]{Example}

\begin{document}

\title{Holomorphic Extension from Weakly Pseudoconcave CR Manifolds}

\author{Andrea Altomani\\
        Research Unity in Mathematics,
        University of Luxembourg\\
        162a, avenue de la Fa\"iencerie,
        L-1511 Luxembourg\\
        E-mail address: andrea.altomani@uni.lu
 \and
        C. Denson Hill\\
        Department of Mathematics,
        Stony Brook University\\
        Stony Brook NY 11794 (USA)\\
        E-mail address: dhilll@math.sunysb.edu
\and
        Mauro Nacinovich\\
        Dipartimento di Matematica,
        II Universit\`a di Roma ``Tor Ver\-ga\-ta''\\
        Via della Ricerca Scientifica,
        00133 Roma (Italy)\\
        E-mail address: nacinovi@mat.uniroma2.it
\and
        Egmont Porten\\
        Department of Mathematics,\
        Mid Sweden University\\
        85170 Sundsvall (Sweden)\\
        E-mail address: Egmont.Porten@miun.se}

\date{}
\maketitle

\begin{abstract}
Let M be a smooth locally embeddable CR manifold, having some CR
dimension $m$ and some CR codimension $d$. We find an improved local
geometric condition on $M$ which guarantees, at a point $p$ on M,
that germs of CR distributions are smooth functions, and have
extensions to germs of holomorphic functions on a full ambient
neighborhood of p. Our condition is a form of weak pseudoconcavity,
closely related to essential pseudoconcavity as introduced in \cite{HN1}.
Applications are made to CR meromorphic functions and mappings.
Explicit examples are given which satisfy our new condition,
but which are not pseudoconcave in the strong sense. These results
demonstrate that for codimension $d > 1$, there are additional phenomena
which are invisible when $d = 1$.
\end{abstract}

\section{Introduction}
The goal of the present article is to give improved geometric conditions on a
generic CR manifold $M\subset\C^n$ which guarantee that all local CR functions
extend holomorphically to a full neighborhood of a given point. This is well known
to be true for strictly pseudoconcave CR manifolds, i.e.~in the case where the
Levi form has one negative eigenvalue in each characteristic conormal direction.
For hypersurfaces this is a classical result of H.~Kneser and H.~Lewy; for $M$ of higher codimension
it was proved independently by a number of authors (see \cite{BP}, \cite{NV}
for CR distributions, and \cite{HN4} for two
different proofs, one very short). Despite numerous efforts, the general problem
to characterize those weakly pseudoconcave manifolds for which one
has extension to full neighborhoods is still far from being completely
understood, even for real analytic hypersurfaces.

Subtle sufficient conditions (sector and ray property) are known for weakly
pseudoconcave {\it hypersurfaces} of finite type (see \cite{BT2}, \cite{FR} for results
and references). In {\it higher codimension} there are several options to approach the
weakly pseudoconcave case. Manifolds which are Levi flat at the reference point
to a certain order, and have all relevant concavity in the generalized Levi form
determined by the next-order terms are studied in \cite{Bo}. Here we aim at the
opposite case where effects of different orders (counted with respect to bracket
length) are combined. Actually our main motivation stems from homogeneous CR manifolds
which biholomorphically look the same near every point. These higher codimensional
homogeneous CR manifolds are abundant, occurring naturally in mathematics, and they
have a strong tendency to be weakly pseudoconcave (see \cite{MN1}, \cite{MN2}, \cite{MN3},
\cite{MN4}, \cite{AMN}). The main result of the present article reveals that there are
additional phenomena which are invisible in codimension one, and indicates that
finite type together with a suitable
notion of weak pseudoconcavity should imply extension to a full
neighborhood. To avoid confusion, we stress that the problem
under consideration is different in nature from the problem of
holomorphic wedge extension, for which a definitive answer is known
(\cite{T1}, \cite{Tu1}, \cite{J}, \cite{M}, see also \cite{MP2}).
In fact, this definitive answer was obtained without having explicit
control on the directions of extension, which is crucial for the
problem at hand.

Let $M\subset\C^n$ be a smooth CR manifold. We denote by $J$ the
complex structure tensor on $T\C^n$, by $HM=TM\cap JTM$ the holomorphic tangent
bundle of $M$, which is the real subbundle of $TM$ invariant under $J$,
and by $H^0 M\subset T^*M$ the characteristic bundle, defined
fiberwise as the annihilator of $H M$. We define the {\it vector valued} Levi form
\[
\cL_{M,p}=\cL_{p}:H_p M\times H_p M\rightarrow \C\otimes (T_p M/H_p M)
\]
by
\[
\cL_{M,p}(X,Y)=-[\tilde{X},J\tilde{Y}](p)+i[\tilde{X},\tilde{Y}](p)
\mbox{ mod }\C\otimes H_p M,
\]
where $\tilde{X},\tilde{Y}\in\Gamma(M,HM)$ are smooth extensions of $X,Y$,
respectively. The {\it usual} Levi form ${\cal L}_{M,\xi}={\cal L}_\xi$ is parameterized by
the characteristic codirections $\xi\in H^0_p M$. It is defined by
$\cL_{M,\xi}(X,Y)=\xi(\cL_{p}(X,Y))$, for $\xi\in H^0_p M$, $X,Y\in H_p M$
(where we read $\xi$ as a form on $\C\otimes (T_p M/H^0_p M)$ in the canonical way).
Most often we shall work with the corresponding real hermitian forms
$\cL_{p}(X)=\cL_{p}(X,X)\in T_p M/H_p M$ and $\cL_{\xi}(X)=\cL_{\xi}(X,X)$.
The reader should take note of the subtle difference in notation between
$\cL_{p}(X)$ and $\cL_{\xi}(X)$; the former is vector valued, and the latter
is scalar valued.

A CR manifold $M$ is {\it strictly (weakly) pseudoconcave} at a point $p\in M$ if
for every $\xi\in H^0_p M$, $\xi\not=0$, $\cL_{\xi}$ has a negative (nonpositive)
eigenvalue. Replacing $\xi$ by $-\xi$, we see that $\cL_{\xi}$ has actually eigenvalues
of both signs (in the strictly pseudoconcave case). Following \cite{HN1}, we call
$M$ {\it trace pseudoconcave} at $p\in M$ if for every $\xi\in H^0_p M$,
$\cL_{\xi}$ is either zero or has eigenvalues of both signs. Trace pseudoconcavity
isolates one of the properties of essential pseudoconcavity introduced in \cite{HN1}.
We refer to this article for background information.

Let ${\cal G}_1$ be the sheaf of germs of smooth (real) CR vector fields on M
(i.e.~sections of $HM$).
For every positive integer $k$ we define inductively ${\cal G}_{k+1}$ as the sheaf
generated by ${\cal G}_k$ and $[{\cal G}_1,{\cal G}_k]$. Let
$G_{k,p}\subset T_p M$ be the vector space generated by pointwise evaluations
of germs in ${\cal G}_k$ at a point $p\in M$. We say that $M$ is of
{\it kind $k$} at $p$ if $G_{k,p}=T_p M$ but $G_{j,p}\subsetneq T_p M$ for $j<k$.
We say that $M$ satisfies the {\it constant rank condition} if the spaces
$G_{k,p}$ have dimension independent of $p$,
i.e.~if they form vector bundles $G_k=\bigcup_{p\in M}G_{k,p}$.

Now we can formulate our main result.

\begin{theorem}\label{main}
Let $M$ be a smooth generic CR manifold in $\C^n$ and $p_0\in M$.
Assume that in a neighborhood of $p_0$, $M$ is trace pseudoconcave,
satisfies the constant rank condition and is of kind 
less or equal to $3$.
Then for every open neighborhood $U$ of $p_0$ in $M$, there is
an open neighborhood $V$ of $p_0$ in $\C^n$ such that every CR
distribution on $U$ is smooth on $M\cap V$ and has
a unique holomorphic extension to $V$.
\end{theorem}

We emphasize that Theorem \ref{main} reveals a phenomenon which
remains invisible in codimension $1$. In fact, in the hypersurface
case its assumptions imply that $M$ is of kind 2, hence strictly
pseudoconcave. We expect the result to extend to arbitrary finite kind.
In \cite{HN1}, the weak identity principle for CR functions (coincidence
on open sets implies coincidence everywhere) was shown for essentially
pseudoconcave CR manifolds. For those CR manifolds covered by the assumptions
of Theorem \ref{main} our result immediately yields the strong identity
principle (coincidence of Taylor coefficients at some point yields
coincidence in a neighborhood).

Theorem \ref{main} is proved in Sections 2-4. The rest of the article
is devoted to applications, extensions and examples. In Section 5 we observe
that CR manifolds as in Theorem \ref{main} enjoy the local extension property E
introduced in \cite{HN2}. The results in \cite{HN2}, \cite{HN3} yield far
reaching global consequences for fields of CR meromorphic functions on
such manifolds. Here CR meromorphic functions are functions which are
locally representable as fractions of CR functions. An alternative approach
to CR meromorphic mappings originates from work of Harvey and Lawson \cite{HL}.
The idea is to require the graph to look like a CR manifold with appropriate
singularities. In general, extension of such mappings is complicated. Based
on \cite{MP2}, we prove in Section 6 that such CR meromorphic functions extend
meromorphically from manifolds with property E to full ambient neighborhoods
and are in particular representable as local quotients. In Section 7 we
present several classes of homogeneous CR manifolds to which all the local and global
results indicated above apply. These were discovered in a much broader context
(see \cite{AMN}). For the reader's comfort, we give a reasonably self-contained
presentation.

{\bf Acknowledgement:}
The second and the fourth author would like to thank the Mathematisches Institut
der Hum\-boldt-Universit\"at zu Berlin, and Professor J\"urgen Leiterer in particular,
for their kind hospitality. Furthermore we thank Anna Siano for information about the
literature.

\section{Preliminaries}

We will use some standard facts about the bundles
$G_1\subset G_2\subset\ldots$: If $G_k=G_{k+1}$ then all $G_j$,
$j\geq k$, are equal (the proof is an application of the
Jacobi identity). This means in particular that $G_k$ is integrable
in the sense of Frobenius. Moreover the map associating to smooth
sections $X\in\Gamma(U,G_1)$, $Y\in\Gamma(U,G_k)$, the section
$[X,Y]\mbox{ mod }G_k\in\Gamma(U,G_{k+1}/G_k)$ is tensorial,
i.e.~$[X,Y](p)\mbox{ mod }G_{k,p}$ depends only on $X(p)$
and $Y(p)$.

Let us now have a closer look at $G_2$. First we note that
independently of concavity $G_{2,p}/H_p M$ is spanned as a
real vector space by the image  $C_p=\{{\cal L}_p(X):X\in H_p M\}$
of the vector valued Levi form. Indeed, $L=\mbox{span}_{\R}C_p$
is contained in $G_{2,p}/H_p M$ by definition. On the other hand,
polarization shows that
$L\otimes\C=\mbox{span}_{\C}\{{\cal L}_p(X,Y):X,Y\in H_p M\}$.
Since the imaginary part of ${\cal L}_p(X,Y)$ is essentially
$[X,Y]$, we obtain $G_{2,p}/H_p M\subset L$. The above is
equivalent to the fact that $G_{2,p}$ is spanned by $H_p M$
and the preimage of $C_p$ under the canonical projection
$T_p M\rightarrow T_p M/H_p M$.

A simple but crucial observation is that trace pseudoconcavity
allows us to replace linear spans by convex hulls.

\begin{lemma}\label{ch}
Assume that $M$ is trace pseudoconcave at $p\in M$. Then $G_{2,p}/H_p M$
is the convex hull of $C_p=\{{\cal L}_p(X):X\in H_p M\}$.
\end{lemma}

{\bf Proof:}
If the lemma fails, there is a nonzero linear functional $\xi$ on
$G_{2,p}/H_p M$ such that $C_p\subset\{\xi\geq 0\}$. We may extend
$\xi$ to an element of $H^0_p M$. Since $G_{2,p}/H_p M$ is the
linear span of $C_p$, there is an $X\in C_p$ with
$\xi(X)>0$. Any $Y\in H_p M$ with ${\cal L}_p(Y)=X$ satisfies ${\cal L}_\xi(Y)>0$.
But this implies that ${\cal L}_\xi$ has some negative eigenvalue,
in contradiction to $C_p\subset\{\xi\geq 0\}$ $\Box$

To keep track of directions of extension, we will use an analogue of the analytic wave
front set, denoted by $WF_u$. It is defined for CR distributions $u$ via the FBI transform in \cite{S},
see also \cite{T2}. For $U$ open in $M$, let $C\!R(U)$ denote the space
of continuous CR functions defined on $U$.
We do not even have to recall the definition of $WF_u$, since the following basic
properties will suffice for our purposes:
\begin{itemize}
\item [{\bf (a)}] Let $u$ be a CR distribution defined on $U\subset M$. Then $WF_u$ is a cone closed in the
    pointed characteristic bundle $H^0 U\backslash o$ ($o$ denoting the zero section).
\item [{\bf (b)}] $WF_u\cap H^0_p M=\emptyset$ holds if and
    only the CR distribution $u$ extends holomorphically to an ambient neighborhood of $p$.
\item [{\bf (c)}] Let $u\in C\!R(U)$. If CR extension from $U$ holds at
    $(p,X)$, $p\in U$, $X\in T_p M\backslash H_p M$, then
    for any $\xi\in WF_u\cap H^0_p M$ we have $\xi(X)\geq 0$.
\end{itemize}
In {\bf (c)} we use the following terminology: We say that {\it CR extension from $U$} holds at
$(p,X)$, $p\in U$, $X\in T_p M\backslash H_p M$, if there is a ${\cal C}^2$-smooth ($\dim M+1$)-dimensional
CR manifold $\tilde{M}$  attached to $U$ along some $U$-neighborhood $U'$
of $p$ such that (i) for a representative of $X$, $JX$ points into $\tilde{M}$ and
(ii) every $u\in C\!R(U)$ has a continuous extensions to $\tilde{M}\cup U'$ which
is CR on $\tilde{M}$. Neglecting the dependence on $U$, we will sometimes call $(p,X)$
or just $X$, a direction of CR extension.
For {\bf (a)}, {\bf (b)}, see \cite{S}, whereas {\bf (c)} is observed in \cite{T2}.

Theorem \ref{main} is a consequence of the following more precise result which does not
require kind $3$.

\begin{theorem}\label{WF}
Let $M$ be a smooth generic CR manifold in $\C^n$.
Assume that on an open set $U\subset M$, $M$ is
trace pseudoconcave and that $G_2$, $G_3$ are bundles.
Then for every continuous CR function $u$ defined on
$U$, we have $WF_u\subset G_3^\perp$.
\end{theorem}

Theorem \ref{WF} will be proved in the next two sections.
It implies Theorem \ref{main} in the following way:
Kind 3 means that $G_3^\perp$ is the zero bundle near $p_0$.
By {\bf (c)} a continuous CR function $u$ extends holomorphically
to an ambient neighborhood of every $p$ contained in some
neighborhood $U'$ of $p_0$ in $M$. By a standard gluing argument
one obtains extension to an ambient neighborhood $V'$ of $U$,
which may a priori depend on $u$. Note that holomorphic
extension in particular shows that $u\in{\cal C}^\infty(U')$.
Now a Baire category argument as in \cite{HN2} yields extension
to a neighborhood of $V$ whose size only depends on $U$.
This proves the theorem for continuous CR functions.

If $u$ is a CR distribution, we may use a method from \cite{BT1}, \cite{T},
to represent it near $p_0$ as $u=\Delta_M^k f$, where $f$ is a continuous CR
function. Here $k$ is a sufficiently large integer, and $\Delta_M$ is a variant
of the Laplace operator which is defined in an ambient neighborhood of $p_0$
and restricts nicely to $M$. If $\tilde{f}$ is a holomorphic extension of $f$,
then the various properties of $\Delta_M$ imply that $\Delta_M^k \tilde{f}$
is the desired extension of $u$. Hence Theorem \ref{main} follows from
Theorem \ref{WF}.

\section{Proof plan for Theorem \ref{WF}}

Here we will prove Theorem \ref{WF} modulo some more technical results
on CR extension which are postponed to the following section.
Pick some $u\in C\!R(U)$.

{\bf Step 1:} $WF_u\subset G_2^\perp$.
This will follow from trace pseudoconcavity. It is a consequence
of the following lemma which holds without constant rank assumptions.

\begin{lemma}
Assume that $M\subset\C^n$ is trace pseudoconcave at $p\in M$.
Then for every continuous CR function $u$ defined near $p$ we
have $WF_u\cap H^0_{p} M\subset G_{2,p}^\perp$.
\end{lemma}

{\bf Proof:}
Let $\xi\in WF_u\cap H^0_{p} M$. We will consider $\xi$ both as a functional
acting on $T_p M$ and $T_p M/H_p M$.
By \cite{Tu2}, every element $X\in C_p$ can be approximated by directions of
CR extensions $X_j\in T_p M/H_p M$. From property {\bf (c)} and continuity we get $\xi(X)\geq 0$.
It follows that $\xi$ is nonnegative on the convex hull of $C_p$. Since
this convex hull is the vector space $G_{2,p}/H_p M$, $\xi$ vanishes on
$G_{2,p}$ $\Box$

Notice that the lemma together with property {\bf (b)} already imply extension to a full
neighborhood for strictly pseudoconcave CR manifolds.

{\bf Step 2:} $WF_u\subset G_3^\perp$.
By Step 1 it suffices to show that a covector $\xi_0$ that annihilates
$G_2$ but not $G_3$ is not contained in $WF_u$. Let $p_0\in U$ denote
the base point to which $\xi_0$ projects. It is our aim to show
$\xi_0\notin WF_u$ by constructing an appropriate CR extension and
applying property {\bf (c)}.

Since the image $C_{p_0}$ of the vector valued Levi form spans $G_{2,p_0}/H_{p_0} M$
(see the remarks before Lemma \ref{ch}),
we may select vectors $X_1,\ldots,X_k\in H_{p_0} M$ such that the vectors
$\tilde{Y}_j=\cL_{M,p_0}(X_j)$, $j=1,\ldots,k$, form a basis of $G_{2,p_0}/H_{p_0} M$.
Extending the $X_j$ smoothly to CR vector fields defined near $p_0$, we
obtain a local basis $\tilde{Y}_j(p)=\cL_{M,p}(X_j)$ of $G_{2}/HM$. Set
$Y_j=[J X_j,X_j]$, and choose a local basis $Z_1,\ldots,Z_{2m}$ of $HM$.
Then the $Z_i$ form together with the the $Y_j$ a local basis of $G_2$.

First we claim that $G_3$ is spanned in some neighborhood
of $p_0$ by the $Z_i$, $Y_j$, together with the brackets $[Z_i,Y_j]$.
Indeed, by definition $G_{3,p}$ is spanned by $G_{2,p}$ and vectors of the
form $[Z,Y](p)$ where $Z\in{\cal G}_{1,p}$, $Y\in{\cal G}_{2,p}$ for $p$
near $p_0$. Around $p$ we may write $Z=\sum z_i Z_i$, $Y=\sum \tilde{z}_i Z_i+\sum y_j Y_j$,
with smooth coefficients $z_i$, $\tilde{z}_i$, $y_j$. This yields
\[
[Z,Y]=\sum z_i\, y_j [Z_i, Y_j]+R,
\]
where $R$ is a germ in ${\cal G}_{2,p}$. This proves the claim.

In the sequel, we will only need the following consequence:
Since $\xi_0$ does not annihilate $G_{3,p_0}$, there are $i_0$, $j_0$ such that
$\xi_0([Z_{i_0},Y_{j_0}](p_0))\not=0$. For notational convenience we will write
from now on $Y=Y_{j_0}$, $Z=Z_{i_0}$.

The following proposition, which will be proved in Section 4, yields CR extension
at $(p_0,Y(p_0))$ realized by a CR manifold to which $Y$ is complex tangent {\it
in a neighborhood} of $p_0$ in $M$.

\begin{proposition}\label{attachmf}
Let $M\subset\C^n$ be a smooth generic CR manifold of CR dimension $m$ and codimension $d$.
Let $p_0\in M$ and let $U$ be an open neighborhood of $p_0$ in $M$. Let $X$ be a smooth CR vector
field on $U$ with $\cL_{M,p_0}(X)\not=0$. Then there is a local ${\cal C}^4$-smooth generic
CR manifold $\tilde{M}$ of dimension $\dim M+1$ with the following properties:

{\bf (a)} $M\cap\tilde{M}$ is a neighborhood of $p_0$ in $U$ and
$\tilde{M}\backslash M$ has two connected components $\tilde{M}^\pm$.\\
{\bf (b)} The distribution on $M\cap\tilde{M}$ spanned by $HM$ and $[JX,X]$
coincides with $(H\tilde{M}\cap TM)|_{M\cap\tilde{M}}$.\\
{\bf (c)} Denoting by $\tilde{M}^+$ the side into which $J[JX,X](p_0)$ points,
we have CR extension from $U$ to $\tilde{M}^+$.
\end{proposition}

Since $\xi_0$ annihilates $H_{p_0}\tilde{M}\cap T_{p_0} M$, there is a unique extension
$\tilde{\xi}_0\in H^0_{p_0}\tilde{M}$. We claim that $\cL_{\tilde{M},\tilde{\xi}_0}$
has eigenvalues of both signs. Indeed, both $Z$ and $Y$ may be extended to
${\cal C}^3$-smooth CR vector fields on $\tilde{M}$ which we denote by the
same symbols. For $Y$ we use here that $Y$ is {\it complex} tangent to $\tilde{M}$
along $M$ by Proposition \ref{attachmf}. The fact that $\xi_0$
annihilates the image of the Levi form of $M$ implies
$\cL_{\tilde{M},\tilde{\xi}_0}(Z(p_0))=\cL_{M,\xi_0}(Z(p_0))=0$.
Since $\tilde{\xi}_0([Z,Y](p_0))\not=0$ is essentially the imaginary
part of the sesquilinear Levi form, $Z$ does not lie in the
kernel of $\cL_{\tilde{M},\tilde{\xi}_0}$. Hence
$\cL_{\tilde{M},\tilde{\xi}_0}(\cdot)$ takes both positive and
negative values in every neighborhood of $Z(p_0)$ in $H_{p_0}\tilde{M}$.
The claim follows.

Let $X^+$, $X^-$ be positive and negative eigenvectors of
$\cL_{\tilde{M},\tilde{\xi}_0}$, respectively.
Since all vectors in $H_{p_0}M$ are null vectors of
$\cL_{\tilde{M},\tilde{\xi}_0}$, we may assume, after replacing $X^\pm$ by
an appropriate complex multiple if necessary, that $X^\pm$ both point into $\tilde{M}^+$.
But now the following proposition, whose proof is postponed to Section 4,
yields for $M$ two directions of CR extension $Y^\pm\in T_{p_0}M/H_{p_0}M$ with $\xi_0(Y^+)>0$
and $\xi_0(Y^-)<0$. This implies $\xi_0\not=WF_u$ by property {\bf (c)}.
Hence the proof of Theorem \ref{WF}, and thus also of Theorem \ref{main},
will be complete as soon as we have shown the following.

\begin{proposition}\label{newdir}
Let $M$, $\tilde{M}$, $\tilde{M}^\pm$ and $p_0\in M$ be as above.
Assume that there is $\tilde{\xi}\in H^0_{p_0}\tilde{M}$
and $X\in H_{p_0}\tilde{M}\backslash T_{p_0}M$ pointing into $\tilde{M}^+$
such that $\cL_{\tilde{M},\tilde{\xi}}(X)>0$. Then there is a direction
$Z\in T_{p_0} M/H_{p_0} M$ of CR extension from $U$ satisfying $\tilde{\xi}(Z)>0$.
\end{proposition}

\section{CR extension}

In this section we will prove Propositions \ref{attachmf} and \ref{newdir}.
First we recall some basic material on Bishop discs introduced in the
seminal paper \cite{B}. Consider a generic CR manifold $M\subset\C^n$
of CR dimension $m$ and codimension $d$. A ${\cal C}^*$-smooth analytic
disc is a mapping $A(\zeta)=(Z(\zeta),W(\zeta))\in{\cal C}^*(\overline{\Delta},\C^n)
\cap{\cal O}(\Delta,\C^n)$ where $\Delta=\{\zeta\in\C:|\zeta|<1\}$.
We say that $A$ is attached to $M$ if $A(S^1)\subset M$.

For $|\mu|\leq 1$ let ${\cal T}_{\mu}$ denote the Hilbert
transform of a function $U:S^1\rightarrow\R$ to its harmonic conjugate
${\cal T}_{\mu} U$, normalized by the condition that the harmonic extension
of ${\cal T}_{\mu} U$ vanishes at $\zeta=\mu$. It is known that
${\cal T}_{\mu}$ is a continuous linear operator on the H\"older
spaces ${\cal C}^{k,\alpha}(S^1,\R)$ if $k\in\N_0$, $0<\alpha<1$. We shall
use the same notation for the Hilbert transform applied componentwise to
vector valued functions of ${\cal C}^{k,\alpha}(S^1,\R^d)$.

We work in coordinates
\begin{equation}\label{v=h}
z_1=x_1+iy_1,\ldots,z_m=x_m+iy_m,w_1=u_1+iv_1,\ldots,w_d=u_d+iv_d
\end{equation}
centered at the origin in which $M$ is locally given as a graph
$v=h(z,u)$ with $h(0)=0$, $dh(0)=0$. The Bishop equation is the
nonlinear functional equation
\begin{equation}\label{bishop}
U=-{\cal T}_\mu(h(Z,U))+u.
\end{equation}
Here $U=U(\zeta)$ is the unknown function mapping the unit circle
$S^1=\{|\zeta|=1\}\subset\C$ to $\R^d$, whereas $Z=Z(\zeta)$ is the
boundary value of a given holomorphic function from $\Delta$ to $\C^m$
and $u\in\R^d$ is a prescribed vector.
It is known that the Bishop equation can be solved in ${\cal C}^{k,\alpha}(S^1,\R^d)$
provided the data $Z(\zeta)$ are ${\cal C}^{k,\alpha}$-small and
$|u|$ is small (see \cite{MP4} for detailed information).
The solution $U$ corresponds to a unique analytic disc $A:\overline{\Delta}\rightarrow\C^n$
whose restriction to $S^1$ is $(Z(\zeta),U+ih(Z(\zeta),U(\zeta)))$. More precisely,
$A$ is holomorphic on $\Delta$ and ${\cal C}^{k,\alpha}$-smooth up to $S^1$. By
construction, $A$ is attached to $M$.

{\bf Proof of Proposition \ref{attachmf}:}
The construction of $\tilde{M}^+$ will be a refinement of \cite[proof of Theorem 9.1]{HT}.
We may assume $p_0=0$ and choose coordinates $(z,w)$
as in (\ref{v=h}). After appropriate rotations and dilations in $z$ and $w$,
we may furthermore assume ${\cal L}_{M,0}(\frac{\partial}{\partial x_1})=
\frac{\partial}{\partial u_1}\mbox{ mod }H_0 M$.
We will construct $\tilde{M}^+$ as the union of real curves $\gamma_p(s)$,
$s\in[0,\epsilon_1)$, having their initial points at $\gamma_p(0)=p\in M$
and such that the segments $\gamma_p((0,\epsilon_1))$ foliate $\tilde{M}^+$.

Let us start by the construction of $\gamma_0$. Following \cite{HT}, we first
simplify the defining equations by removing some of the pure terms. After a
change of the $w_1$-variable
\[
w_1\mapsto w_1-i\Big(
\sum_{j,k=1}^m\frac{\partial^2 h_1}{\partial z_j \partial z_k}(0) z_j z_k
+\sum_{j=1}^d\frac{\partial^2 h_1}{\partial z_1\partial u_j}(0)z_1 w_j
+\sum_{j,k=1}^d\frac{\partial^2 h_1}{\partial u_j \partial u_k}(0) w_j w_k
\Big),
\]
we have
\[
h_1(z,u)=\sum_{j,k=1}^{m} a_{j,k} z_j \overline{z}_k + O_3(z,u),
\]
with $a_{1,1}=1$. In particular, we still have
${\cal L}_{M,0}(\frac{\partial}{\partial x_1})=\frac{\partial}{\partial u_1}$.

For a small $\epsilon_1>0$, to be specified later, and a parameter $0\leq t<\sqrt{\epsilon_1}$,
we let $U_t$ be the solutions of the parameter-dependent Bishop equation
\begin{equation}\label{bishop}
U_t=-{\cal T}_0(h(W_t,U_t)),\mbox{ where }W_t=(t\zeta,0\ldots,0).
\end{equation}
Since the data are smooth, we can solve this in the H\"older spaces
${\cal C}^{k,\alpha}(S^1,\R^d)$ ($k\geq 1$, $0<\alpha<1$) for $\epsilon_1$
sufficiently small. Let the $A_t(\zeta)$ be the corresponding holomorphic disc
and set $\gamma_0(s)=A_{\sqrt{s}}(0)$. One reads off from (\ref{bishop}) that the curve
$\gamma_0$ starts at the origin and runs in the space $i\R_v^d =\{z=u_1=\ldots=u_n=0\}$.
It is shown in \cite{HT} that (i) $\frac{d\gamma_0}{ds}(0)$ is a positive multiple of
$\frac{\partial}{\partial v_1}$ and that (ii) $\gamma_0$ is on $[0,\epsilon_1)$ as smooth
as we please if $\epsilon_1$ is sufficiently small.

The curves $\gamma_p$ will be obtained by varying the initial point in this construction.
To this end we first produce a family of coordinates $(z_p,w_p)$ centered at $p$
such that $T_{p=0}M=\{v_p=0\}$. Clearly this can be achieved by an affine linear change of
coordinates $\Phi_p$ which depends smoothly on $p$, where $p$ ranges in some 
$M$-neighborhood of the origin.
Thus $(z_p,w_p)$ are holomorphic for $p$ fixed, but are only smooth in $p$. Next we rotate
and dilate in such a way that $d\Phi_p(X)$ transforms to $\frac{\partial}{\partial x_1}$ at
the origin and such that ${\cal L}_{M,0}(\frac{\partial}{\partial x_1})=\frac{\partial}{\partial u_1}$.
Then we modify $(z_p,w_p)$ again as before in order to simplify second-order
terms. Clearly all this can be done by a family of local biholomorphisms depending smoothly
on the parameter $p$. In these last coordinates, which we still denote by $(z_p,w_p)$,
we construct a curve $\tilde{\gamma}_p(s)$, $0\leq s<\epsilon_1$,
verbatim by the same construction as above. Thus $\tilde{\gamma}_p(s)$ starts at the
origin in $(z_p,w_p)$-space, which corresponds to the point $p$, and its time derivative
at $s=0$ is a positive multiple of $\frac{\partial}{\partial v_1}$.
The desired curve $\gamma_p$ with $\gamma_p(0)=p$ is then obtained
by reversing the coordinate transformations. The nature of the process
implies that $\frac{d\gamma_p}{ds}(0)$ is a positive multiple of $X(p)$.

We claim that $\tilde{M}^+\cup W=\bigcup_{p\in W,\,0\leq s<\epsilon_1}\gamma_p(s)$
is a manifold with boundary $W$ attached to $M$ along $W$,
provided $W$ is a sufficiently small open neighborhood
of the origin in $M$, and $\epsilon_1$ is sufficiently small. Actually, optimal
regularity results for the Bishop equation yield that the solution, which depends on
a finite dimensional set of parameters, has arbitrarily small loss of smoothness with respect
to the parameters (\cite{Tu2}, see also \cite{MP4}). Thus we have that data depending
${\cal C}^{k,\alpha}$-smoothly
on all variables and parameters lead to ${\cal C}^{k,\alpha-0}$-smooth solutions.
As the smoothness in $s$ is as good as needed and we can assume $k$ as large as
we please, the mapping $(p,s)\mapsto\gamma_p(s)$ is as smooth as we please.
First, the inverse function theorem implies the claim about $\tilde{M}^+\cup W$.
Second, we may extend $(p,s)\mapsto\gamma_p(s)$ to $W\times(-\eta,\epsilon_1)$ for some
$0<\eta<<1$, and obtain an extended manifold
$\tilde{M}=\bigcup_{p\in W,-\eta<s<\epsilon_1}\gamma_p(s)$ as required in {\bf (a)}.
The before mentioned fact that $\frac{d\gamma_p}{ds}(0)$ is a positive
multiple of $X(p)$ completes the proof of {\bf (b)}.

The proof of {\bf (c)} is standard. One chooses an $M$-neighborhood $W'\subset\subset U$
of $p_0$ so small that every $u\in C\!R(U)$ can be uniformly approximated
on $\overline{W'}$ by holomorphic polynomials $P_j$, using the Baouendi-Treves
approximation theorem \cite{BT1}. Restricting domains again, one restricts the
above construction using only discs attached to $W'$. Then the
maximum modulus principle applied to the discs implies that the $P_j$ converge uniformly
on $\tilde{M}^+\cup W$ to a continuous function $u_+$, which is CR on $\tilde{M}^+$.
The proof of Proposition \ref{attachmf} is complete $\Box$

{\bf Proof of Proposition \ref{newdir}:}
Let us first sketch the geometrical idea. Since
$Y={\cal L}_{\tilde{M},p_0}(X)$ is not complex tangent to $\tilde{M}$,
it is classical (see \cite{HT}) that appropriate ${\cal C}^{2,\alpha}$-small discs
attached near $p_0$ to $\tilde{M}$ whose $z$-coordinates are parallel to $\C X$ are
nontangent to $\tilde{M}$ along their boundaries, sticking out along directions
which are approximately $Y$ modulo $H_{p_0}\tilde{M}$. This would be enough
in order to extend CR functions defined on all of $\tilde{M}$. In the case
at hand, we are only allowed to use discs attached to $\tilde{M}^+\cup W$.
We will construct a family of discs whose boundaries touch $M$ quadratically
in exactly one point and obtain CR extension to a $(\dim M+1)$-dimensional manifold
(distinct from $\tilde{M}^+$)
contained in the union of the discs. Since essentially the same construction is explained
in great detail in \cite[Section 5]{MP3}, it will suffice to give a concise review
of what has to be done.

We will first construct a single disc attached to $\tilde{M}^+\cup W$ whose boundary
touches $M$ at $p_0$. Choose local coordinates
\[
z_1,\ldots,z_{m+1},w_1=u_1+iv_1,\ldots,w_{d-1}=u_{d-1}+iv_{d-1}
\]
centered at $p_0$ such that $\tilde{M}$ is locally given as a graph
$v=\tilde{h}(z,u)$, with $\tilde{h}(0)=0$, $d\tilde{h}(0)=0$.
After convenient rotations, we may assume that
$X=\frac{\partial}{\partial x_1}$,
${\cal L}_{\tilde{M},0}{X}=\frac{\partial}{\partial u_1}\mbox{ mod }H_0\tilde{M}$
and $T_0 M=i\R_{y_1}\oplus\C^{m}_{z_2,\ldots,z_{m+1}}\oplus
\R^{d-1}_{u_1,\ldots,u_{d-1}}$ (by multiplying by some $\zeta\in S^1$,
we have rotated $X$ so that $JX\in T_0 M)$). For technical reasons, we also
arrange that
\begin{equation}\label{hN}
\tilde{h}(z_1,0,\ldots,0,u_1\ldots,u_{d-1})=c|z_1|^2+O(|(z,u)|^3),\,\,\, c>0,
\end{equation}
by eliminating pure terms of second order.

For $Z_r(\zeta)=(r(1-\zeta),0,\ldots,0)$, $0<r<<1$, let
\[
A_r(\zeta)=(Z_r(\zeta),U_r(\zeta)+iV_r(\zeta))
=(Z_r(\zeta),U_r(\zeta)+ih(Z_r(\zeta),U_r(\zeta)))
\]
be the analytic disc obtained by solving the Bishop equation
\begin{equation}\label{bishop2}
U_r=-{\cal T}_1 (h(Z_r,U_r)).
\end{equation}
Note that $A_r(1)=0$. The crucial point is that $r(1-\zeta)$ lies in the
right $z_1$-halfplane and touches the imaginary axis quadratically at the origin,
and that the curvature of its boundary at the origin becomes large for
$r\rightarrow 0$.
First, it is proved in \cite[5.5]{MP3} that for $r>0$ small, we have
$A_r(\overline{\Delta}\backslash\{1\})\subset \tilde{M}^+$ and that
$A_r(S^1)$ touches $M$ quadratically at $0$. There is a geometric
estimate for admissible $r$ which is stable under ${\cal C}^{2,\alpha}$-small
deformations of $\tilde{h}$.
Second, the usual classical argument (see \cite{HT}) based on (\ref{hN}) shows that
\begin{equation}\label{discdirection}
\left. -\frac{\partial V_r}{\partial\lambda}\right|_{\zeta=1}=r\,(c',0,\ldots,0)+o(r),
\end{equation}
with some $c'>0$ independent of $r$ (here $\zeta=\lambda+i\eta)$.

Let us take for granted, for the moment, that we can construct for a small fixed
$r$ a ($\dim M+1$)-dimensional manifold $\hat{M}^+$ attached to $M$, and containing
the image $A_r((1-\epsilon,1))$ of the segment $(1-\epsilon,1)$, such that CR functions
extend from $M$ to $\hat{M}^+$. Then (\ref{discdirection}) yields for the normalized
outgoing direction at the origin
\[
-\frac{\partial v_r}{\partial\lambda}(1)
\Big/\Big|\frac{\partial v_r}{\partial\lambda}(1)\Big|
\rightarrow (c',0,\ldots,0),\mbox{ as }r\rightarrow 0.
\]
But this implies that we can approximate the element defined by
$\frac{\partial}{\partial v_1}$in $T_0 M/H_0 M$ by directions
of CR extensions. Since $\tilde{\xi}(\frac{\partial}{\partial v_1})>0$
we have found a direction of extension $Z$ as desired. Now since a sufficiently
small $r>0$ has been fixed, we drop it from the notation and write $A$.

It remains {\bf (i)} to construct the manifold $\hat{M}^+$ and {\bf (ii)} to
establish CR extension from $M$ to $\hat{M}^+$.
The method for getting {\bf (i)} is very similar to that of the proof of Proposition
\ref{attachmf}. We construct a family of local holomorphic coordinates
\[
z_{1,p},\ldots,z_{m+1,p},w_{1,p},\ldots,w_{d-1,p}
\] 
coinciding with the above coordinates for $p=0$ and satisfying the following properties:
\begin{itemize}
\item [\bf a)] $(z_p,w_p)$ depend ${\cal C}^3$-smoothly on the parameter $p$,
which ranges in a small $M$-neighborhood $W''$ of the origin.
\item [\bf b)] For $p$ fixed, $(z_p,w_p)$ are holomorphic and centered at $p$,
and we have $T_0 M=i\R_{y_{1,p}}\oplus\C^{m}_{z_{2,p},\ldots,z_{m+1,p}}\oplus
\R^{d-1}_{u_{1,p},\ldots,u_{d-1,p}}$.
\end{itemize}
We apply the above construction with dependence on the parameter $p$. This
yields discs $A_p$ attached to $\tilde{M}^+\cup W'$ and touching $M$ quadratically
at $A_p(1)=p$. Using regularity results for the Bishop equation (\cite{Tu2}, also
\cite{MP4}), we see that $A_p(\zeta)$ depends ${\cal C}^{2,1/2}$-smoothly
on $p$ and $\zeta\in\overline{\Delta}$.
After a further shrinking of $W''$ and with a smaller $\epsilon$, we obtain a
manifold $\hat{M}^+\cup W''=\bigcup_{p\in W''}A_p((1-\epsilon,1])$ as required in {\bf (i)}.

To show {\bf (ii)} we just have to remember from the construction of
$\tilde{M}^+$ that every $u\in C\!R(U)$ can be approximated by holomorphic
polynomials $P_j$ which converge uniformly on $\overline{W'}\cup\tilde{M}^+$.
Hence the $P_j$ also converge on $\hat{M}^+\cup W''$ to a continuous CR function
extending $u$. This yields {\bf (ii)} and completes the proofs of Proposition \ref{newdir},
and also of Theorems \ref{main} and \ref{WF} $\Box$

\section{Applications to CR meromorphic functions}

One of our motivations is to find concrete applications for the Siegel-type
theorems proved in \cite{HN2}, \cite{HN3}. In these papers, far reaching
global consequences for the field of CR meromorphic functions are proved
for CR manifolds satisfying a local extension property E. For a ${\cal C}^\infty$-smooth
generic CR submanifold $M$ of a complex manifold $N$, property E means that the canonical
restriction mapping ${\cal O}_{N,p}\rightarrow{\cal CR}_{M,p}$ is
surjective for every $p\in M$. Here ${\cal O}_{N}$ and ${\cal CR}_{M}$
denote the sheaf of germs of holomorphic functions on $N$ and the sheaf
of germs of ${\cal C}^\infty$-smooth CR functions on M, respectively, and
${\cal O}_{N,p}$, ${\cal CR}_{M,p}$ are their stalks at $p$.
We obtain immediately that a manifold satisfying the assumptions of
Theorem \ref{main} in a coordinate neighborhood of each of its points 
has property E. Hence we
obtain all results proved in \cite{HN2}, \cite{HN3} by carrying the
local situation studied here to general manifolds.

Actually the main results in \cite{HN2}, \cite{HN3} concern CR meromorphic
functions rather than CR functions. Similarly as ordinary meromorphic functions,
we define {\it CR meromorphic functions in the usual sense} on $U\subset M$ as functions
which are defined on a dense open subset of $U$ and can be represented near {\it every} point
$p\in U$ as the quotient $p/q$ of ${\cal C}^\infty$-smooth CR functions $p$, $q$,
where $q$ does not vanish identically on any nonempty open subset. If $M$ has
property E, every CR meromorphic function on $U$ is the restriction of a
meromorphic function defined on some ambient neighborhood of $U$ in $N$.

Let $M$ be a smooth compact locally embeddable CR manifold of CR dimension $m$
and CR codimension $d$, which at each point satisfies the hypotheses of Theorem \ref{main}.
Then the field ${\cal K}(M)$ of CR meromorphic functions on $M$ has transcendence
degree $k\leq m+d$. If $f_1,\ldots,f_k$ is a maximal set of algebraically independent
CR meromorphic functions on $M$, then ${\cal K}(M)$ is a simple finite algebraic
extension of the field $\C(f_1,\ldots,f_k)$ of rational functions of the
$f_1, f_2,\ldots,f_k$. Assuming that $M$ is connected, there
is also an equivalence between the algebraic dependence over $\C$, and the analytic
dependence, of a finite set of CR meromorphic functions in ${\cal K}(M)$. When
$M$ has a projective embedding there is an analogue of Chow's theorem, and ${\cal K}(M)$
is isomorphic to the field ${\cal R}(Y)$ of rational functions on an irreducible
projective algebraic variety $Y$, and $M$ has a CR embedding in $\mbox{reg }Y$.
For details, and further applications and remarks, see \cite{HN2}, \cite{HN3}.

\section{CR meromorphic mappings according to Harvey and Lawson}

An alternative notion of CR meromorphic functions and mappings was suggested by
Harvey and Lawson in the context of the complex Plateau problem and studied in
\cite{HL}, \cite{DH}, \cite{DS}, \cite{MP1}, \cite{MP2}. The following
definitions appear in \cite{HL} for hypersurfaces and in \cite{DS} for CR manifolds
of arbitrary codimension. Let $M$ be a smooth generic CR submanifold of a complex manifold $N$ 
of CR dimension $m$ and codimension $d$, and let $X$ be an arbitrary complex manifold.
Then a {\it CR meromorphic mapping $F$ in the sense of Harvey-Lawson} of
an open $U\subset M$ with values in $X$ is given by a triple $(F,{\cal D}_F,\Gamma_F)$
with the following properties:

\begin{itemize}
  \item[{\bf (a)}] ${\cal D}_F$ is an open dense subset of $U$,
  \item[{\bf (b)}] $F:{\cal D}_F\rightarrow X$ is a ${\cal C}^1$-smooth CR mapping,
  \item[{\bf (c)}] the closure of the graph of $F$ in $U\times X$ equals $\Gamma_F$
  and is a local scarred CR cycle of CR dimension $m$ and dimension $2m+d$ in $N\times X$.
\end{itemize}
In {\bf (c)} we mean that $\Gamma_F$ is of locally finite $(\dim M)$-dimensional Hausdorff measure
and contains a closed subset $\sigma$ (the scar set) of $(\dim M)$-dimensional Hausdorff measure zero
such that (i) $\Gamma_F\backslash\sigma$ is a ${\cal C}^1$-smooth CR manifold of same dimension
and CR dimension as $M$ and (ii) in a neighborhood of every $(p,x)\in\Gamma_F$,
integration over $\Gamma_F\backslash\sigma$ yields a {\it closed} current
(see \cite{HL}, \cite{DS}, \cite{MP1} for full details).
If $X$ equals $\P^1$, the complex projective line, we also speak of CR meromorphic functions
in the sense of Harvey-Lawson.

Meromorphic extension of these CR meromorphic mappings is technically complicated because
of the presence of the scar set. Actually it requires a certain machinery to derive the
counterpart of Theorem \ref{main} for CR meromorphic functions in the sense of Harvey-Lawson.
We can prove this for all manifolds with property E.

\begin{theorem}\label{mero}
Let $M$ be a smooth generic CR submanifold of a complex manifold $N$ with property E. 
Then every CR meromorphic function $F$ in the sense of Harvey-Lawson, 
defined on an open set $U\subset M$ and with values in $\P^1$
has a meromorphic extension $\tilde{F}$ to a neighborhood $V$ of $U$ in $N$. More precisely,
the graph $\Gamma_{\tilde{F}}\subset V\times\P^1$ of $\tilde{F}$ satisfies
$\Gamma_{\tilde{F}}\cap(U\times\P^1)=\Gamma_F$. In particular, $F$ is a CR meromorphic function
in the usual sense.
\end{theorem}

A CR manifold $M$ is called {\it minimal at a point $p\in M$} (in the sense of Tumanov)
if there is no germ of a CR manifold $N\subset M$ of the same CR dimension as $M$ and of
lower dimension than $M$ containing $p$.

\begin{lemma}\label{minimal}
Let $M$ be a smooth generic CR manifold in $\C^n$ with property E.
Then $M$ is minimal at every point $p\in M$.
\end{lemma}

{\bf Proof:}
Let us assume that $M$ is not minimal at $p_0$. By \cite{BR} there is a
smooth CR function $u$ defined on an open neighborhood $U'\subset U$ of $p_0$ which
does not extend holomorphically to any open wedge attached to $M$ near $p_0$, in
contradiction to property E $\Box$

{\bf Proof of Theorem \ref{mero}:}
Let $F$ be a $\P^1$-valued CR meromorphic function in the sense of Harvey-Lawson
defined on $U\subset M$. First we construct a local extension of $F$ to an ambient neighborhood
of a given $p_0\in U$. Because of Lemma \ref{minimal}, \cite[Theorem 1.2]{MP2} gives meromorphic
extension to an open wedge ${\cal W}$ attached to a neighborhood $U'$ of $p_0$ in $M$. More
precisely, there is an open truncated cone $C\subset\C^n$ with vertex at the origin and
a meromorphic function $\tilde{F}$ on ${\cal W}=U'+C$ which attains $F$ as continuous
boundary value on ${\cal D}_F$.

>From property E and a Baire category argument (see \cite{HN2}), it follows that smooth CR functions on $U'$
extend holomorphically to a {\it uniform} ambient neighborhood $V'$ of $U'$. Since the envelope
of meromorphy coincides with the envelope of holomorphy for domains in $\C^n$, functions which
are meromorphic on an arbitrarily thin neighborhood of $U'$ extend meromorphically to $V'$.
For fixed $c\in C$, the rigid translates $U'_\epsilon=U'+\{\epsilon c\}$ approach $U'$ for $\epsilon\downarrow 0$.
To obtain the desired extension to a neighborhood of $p_0$, it suffices to choose $\epsilon\in(0,1)$
so small that $p_0\in V'+\{\epsilon c\}$ and to extend $\tilde{F}$ to $V'+\{\epsilon c\}$.

Now a standard gluing argument yields a meromorphic function $\tilde{F}$ which is defined
on a neighborhood $V$ of $U$ in $\C^n$ and coincides with $F$ on ${\cal D}_F$.
It remains to prove that $\Gamma_{\tilde{F}}\cap (U\times\P^1)=\Gamma_F$. Near points in ${\cal D}_F$
this is obvious. From {\bf (b)} in the definition of CR meromorphic mappings and
the corresponding (well known) property of meromorphic functions, we deduce
$\Gamma_F\subset\Gamma_{\tilde{F}}\cap (U\times\P^1)$. Assume that there is
$(p_0,\zeta_0)\in\Gamma_{\tilde{F}}\cap (U\times\P^1)\backslash\Gamma_F$. This is obviously
impossible if $p_0$ is a point near which $\tilde{F}$ is a smooth mapping. Hence it remains
to consider the case in which $p_0$ lies in the indeterminacy set
$\Sigma_{\tilde{F}}=\{p:\{p\}\times\P^1\subset\Gamma_{\tilde{F}}\}$.

Following \cite{DS}, we also consider the indeterminacy set $\Sigma_F$ of $F$ defined by
\[
\Sigma_F=\{p\in U:\{p\}\times\P^1\subset\Gamma_F\}.
\]
Our assumption on $p_0$ means that $p_0\notin\Sigma_F$.
Let $\alpha$ be a biholomorphism of $\P^1$ mapping $\zeta_0$ to $\infty$.
It is observed in \cite{DS} that the set-valued function $\alpha\circ F$
naturally induces a CR distribution $g$ of order one on a neighborhood
$U_{p_0}$ of $p_0$ in $M$. As $M$ is minimal
in $p_0$, Tumanov's theorem \cite{Tu1} and the usual extension techniques for
CR distributions yield holomorphic extension to an open wedge attached to $M$
at $p_0$ (assuming $g$ as boundary value in the weak sense). Then an argument
with approach manifolds as above yields a holomorphic extension $\tilde{g}$ to
a full neighborhood of $p_0$. Observing that $\alpha\circ F$ and $\tilde{g}$
are smooth and coincide at points of ${\cal D}_F$, we obtain that the equality
$\tilde{F}=\alpha^{-1}\circ\tilde{g}$ holds near $p_0$. In particular, $\tilde{F}$
is smooth near $p_0$, in contradiction to $p_0\in\Sigma_{\tilde{F}}$.
The proof of Theorem \ref{mero} is complete $\Box$

\begin{remark}\rm
It requires only little extra work to derive a corresponding result for
CR meromorphic mappings with values in a projective manifold $X$. Note that
in the general case we can only expect
$\Gamma_F\subset\Gamma_{\tilde{F}}\cap(U\times X)$.
\end{remark}

\section{Homogeneous examples}

At first glance, it may seem hard to find examples of CR manifolds of kind 3
satisfying the conditions of Theorem \ref{main}. However, the theory
of homogeneous CR manifolds provides many of them in a very natural
way (see \cite{MN1}, \cite{MN2}, \cite{MN3}, \cite{MN4}, \cite{AMN}).
We give a concise description of a class of such homogeneous CR manifolds,
and refer to \cite{AMN} for more details.

Let $G^{\mathbb{C}}$ be a complex connected semisimple Lie group,
with Lie algebra $\mathfrak{g}^{\mathbb{C}}$, and $G$ a connected
real form of $G^{\mathbb{C}}$, with Lie algebra $\mathfrak{g}$.
Fix a Cartan subgroup $H$ of $G$, which is maximally noncompact, that
is a Cartan subgroup such that a maximal compact torus in $H$ has
minimal dimension, and denote by $\mathfrak{h}$ and
$\mathfrak{h}^{\mathbb{C}}$ the Lie algebras of $H$ and of its
complexification $H^{\mathbb{C}}$. In the set of roots $\mathcal R=
\mathcal R(\mathfrak{g}^{\mathbb{C}},\mathfrak{h}^{\mathbb{C}})$ choose
a subset $\mathcal R^+$ of positive roots, adapted to $\mathfrak g$
(cf.\ \cite{AMN}, Proposition 6.1), and let $\mathcal B$ be the corresponding set
of positive simple roots. To any subset $\Phi\subset\mathcal B$ we
associate the parabolic subalgebra and subgroup
\begin{align*}
\mathfrak{q}_{\Phi}&=\mathfrak{h}^{\mathbb{C}}
+\sum_{\alpha\in\mathcal R^+}\mathfrak{g}^{\mathbb{C}}_{\alpha}
+\sum_{\substack{\alpha\in{-\mathcal R^+}\\(\mathrm{supp}\,\alpha)\cap\Phi=\emptyset}}
\mathfrak{g}^{\mathbb{C}}_{\alpha},\\
Q_\Phi&=\mathrm{Norm}_{G^{\mathbb{C}}}(\mathfrak{q}_{\Phi}).
\end{align*}
Here $\mathfrak{g}^{\mathbb{C}}_{\alpha}$ is the eigenspace in $\mathfrak{g}^{\mathbb{C}}$
of a root $\alpha$, and $\mathrm{supp}\,\alpha$ is the support of a root $\alpha$ in $\mathcal B$. The group
$Q_\Phi$ has Lie algebra $\mathfrak{q}_{\Phi}$.

The group $G^{\mathbb{C}}$ acts via the adjoint representation on $\mathfrak{g}^{\mathbb{C}}$ and on linear
subspaces of $\mathfrak{g}^{\mathbb{C}}$. Fix a subset $\Phi\subset\mathcal B$ and let $d$ be the dimension
of $\mathfrak{q}_{\Phi}$. The orbit
\[
Y=\{\mathrm{Ad}_{\mathfrak{g}^{\mathbb{C}}}(g)(\mathfrak{q}_{\Phi})\mid g\in G^{\mathbb{C}}\}
\subset \mathrm{Gr}_{d}(\mathfrak{g}^{\mathbb{C}})
\]
through $\mathfrak{q}_{\Phi}$ in the Grassmannian of $d$-planes in $\mathfrak{g}^{\mathbb{C}}$ is the
flag manifold of parabolic subalgebras of $\mathfrak{g}^{\mathbb{C}}$ conjugate by an inner automorphism
to $\mathfrak{q}_{\Phi}$. It is a smooth irreducible projective subvariety of
$\mathrm{Gr}_{d}(\mathfrak{g}^{\mathbb{C}})$, isomorphic to $G^{\mathbb{C}}/Q_\Phi$.
The orbit
\[
M=\{\mathrm{Ad}_{\mathfrak{g}^{\mathbb{C}}}(g)(\mathfrak{q}_{\Phi})\mid g\in G\}
\subset \mathrm{Gr}_{d}(\mathfrak{g}^{\mathbb{C}})
\]
of $G$ through $\mathfrak{q}_{\Phi}$ is a smooth generic CR submanifold of $Y$ which, by our choice of
the Cartan subgroup and of the system of positive roots, is compact, and is called the minimal orbit
of $G$ in $Y$.

If all local CR functions near a point $p$ of $M$ extend to a full neighborhood of $p$ in $Y$, then the
pair $(M,Y)$ has property E of \cite{HN2}, {\it hence the field $\mathcal{K}(M)$ of CR meromorphic
functions is isomorphic to the field $\mathcal{R}(Y)$ of rational functions on $Y$}, because $Y$
is the smallest projective variety containing $M$ (see \cite{HN2}).

\begin{example}\rm
Identify $\mathbb{C}^{6}$, with the standard basis
$\{\mathrm e_j\}_{1\leq j\leq 6}$, with the quaternionic vector space
$\mathbb{H}^3$ by setting, for $\lambda\in\mathbb{C}$,
\begin{equation*}
  \mathbf j\lambda\mathrm e_{2j-1}=\bar\lambda\mathrm e_{2j},\qquad
  \mathbf j\lambda\mathrm e_{2j}=-\bar\lambda\mathrm e_{2j-1}.
\end{equation*}
Consider the complex flag manifold
\[
Y=\left\{\ell_{1}\subset\ell_{3}\subset\ell_5\subset\mathbb{C}^{6}\mid
\dim\ell_{2j-1}=2j-1,\; 1\leq j\leq 3
\right\}.
\]
Then $Y$ is a compact $13$-dimensional complex manifold, homogeneous for the
action of $G^{\mathbb{C}}=\mathrm{SL}(6,\mathbb{C})$. Near the point
\[
o=\big(\langle \mathrm e_1\rangle_{\mathbb{C}},
\langle \mathrm e_1,\mathrm e_{2},\mathrm e_{3}\rangle_{\mathbb{C}},
\langle\mathrm e_1,\mathrm e_{2},\mathrm e_{3},
  \mathrm e_{4},\mathrm e_{5}\rangle_{\mathbb{C}}\big)
\]
it admits a holomorphic chart given by the nonconstant entries of the
matrix
\[
A=\left(\begin{matrix}
  1&0&0&0&0 \\
  z_{1}&1&0&0&0 \\
  z_{2}&0&1&0&0 \\
  z_{3}&z_{6}&z_{9}&1&0 \\
  z_{4}&z_{7}&z_{10}&0&1 \\
  z_{5}&z_{8}&z_{11}&z_{12}&z_{13} \\
\end{matrix}\right).
\]
Let $M$ be the real submanifold of Y
given by
\[
M=\{(\ell_{1},\ell_{3},\ell_{5})\in Y\mid \mathbb{H}\ell_{1}\subset\ell_{3},\;\mathbb{H}\ell_{3}\subset\ell_{5})\}.\]
Then $M$ is homogeneous for the action of the real form
$G=\mathrm{SL}(3,\mathbb{H})$ of $G$, and it is a compact real-analytic
homogeneous generic CR submanifold of $Y$. Denoting by $A_j$ the
$j$-th column of $A$, near the point $o$ the manifold $M$ is defined
by the system of equations
\[
\begin{cases}
  \mathrm{rk}(A_1,A_2,A_3,
  \mathbf{j}A_1)=3,\\
  \mathrm{rk}(A_1,A_2,A_3,A_4,A_5,
  \mathbf{j}A_1,\mathbf{j}A_2,\mathbf{j}A_3)=5,
\end{cases}
\]
which, in the coordinates $\{z_j\}_{1\leq{j}\leq{13}}$, are
\begin{equation*}
  \begin{cases}
    z_{2}-\bar{z}_{6}+z_{1}\bar{z}_{3}+z_{3}\bar{z}_{9}
      -z_{1}\bar{z}_{1}\bar{z}_{6}-z_{1}\bar{z}_{2}\bar{z}_{9}=0,\\
    z_{4}-\bar{z}_{8}+z_{1}\bar{z}_{5}+z_{3}\bar{z}_{11}
      -z_{1}\bar{z}_{1}\bar{z}_{8}-z_{1}\bar{z}_{2}\bar{z}_{10}=0,\\
    z_{10}-\bar{z}_{12}+z_{9}\bar{z}_{11}+z_{11}\bar{z}_{13}
      -z_{9}\bar{z}_{9}\bar{z}_{12}-z_{9}\bar{z}_{10}\bar{z}_{13}=0,\\
    z_{5}+\bar{z}_{7}-z_{1}\bar{z}_{4}-z_{3}\bar{z}_{10}
      +z_{1}\bar{z}_{1}\bar{z}_{7}+z_{1}\bar{z}_{2}\bar{z}_{10}=0,\\
  \end{cases}
\end{equation*}
thus $M$ has CR dimension $5$ and codimension $8$. The holomorphic
tangent space $T^{1,0}_oM$ admits the basis
$\{\partial/\partial{z}_j\}_{j=1,3,9,11,13}$.  The space of Levi forms
has real dimension $6$, and in the basis
$\{\partial/\partial{z}_j\}_{j=1,3,9,11,13}$ it is the space of
Hermitian symmetric matrices of the form
\[
\begin{pmatrix}
  0&\alpha&0&0&0\\
  \bar\alpha&0&\alpha&\beta&0\\
  0&\bar\alpha&0&\gamma&0\\
  0&\bar\beta&\bar\gamma&0&\gamma\\
  0&0&0&\bar\gamma&0
\end{pmatrix},\qquad \alpha,\beta,\gamma\in\mathbb{C},
\]
as the fourth equation does not contribute to the Levi form. In some
characteristic codirections the Levi forms are zero, but in the
remaining characteristic codirections, the Levi forms all
have signature $(1,1)$ or signature $(2,2)$,
and $M$ is trace pseudoconcave.  Moreover $M$ has kind~$3$, indeed
$G_2$ has real dimension $16$ and $G_3$ is the whole tangent space of
$M$.
By Theorem~\ref{main} germs of CR distributions
on $M$ are real-analytic and extend to germs of holomorphic functions on $Y$.

With the notation of \cite{AMN}, $M$ is
the minimal orbit of the simple Lie group of type $\mathrm{A\,II}_5$
associated to the parabolic subalgebra
$\mathfrak{q}=\mathfrak{q}_{\{\alpha_1,\alpha_3,\alpha_5\}}$.
\end{example}

\begin{example}\rm
  The example above can be generalized to the following pairs of
  complex flag manifolds, homogeneous for the action of
  $\mathrm{SL}(2n,\mathbb{C})$, and compact generic CR submanifolds,
  homogeneous for the action of $\mathrm{SL}(n,\mathbb{H})$:
  \begin{equation*}
    \begin{cases}
      Y=\big\{\ell_{2j_1-1}\subset\ell_{2j_2-1}\subset\ell_{2j_3-1}
      \subset\mathbb{C}^{2n}\mid
      \dim\ell_{2j_k-1}=2j_k-1, 1\leq k\leq3\big\}, \\[4pt]
      M=\left\{(\ell_{2j_1-1},\ell_{2j_2-1},\ell_{2j_3-1})\in Y\left|
          \begin{aligned}
            &\dim(\ell_{2j_k-1}\cap\mathbb H\ell_{2j_k-1})=2j_k-2,
            &&1\leq k\leq3\\
            &\;\mathbb H\ell_{2j_k-1}\subset\ell_{2j_{k+1}-1},&&1\leq
            k\leq2
          \end{aligned}
        \right.\right\},\\[4pt]
      \text{for $1=j_1<j_2<j_3\leq n$.} 
    \end{cases}
  \end{equation*}
  In this case the complex dimension of $Y$, the CR dimension of
  $M$, and the CR codimension of $M$ are:
  \begin{align*}
    \dim_{\mathbb{C}}Y&=(2n-1)
    +2(j_2-1)(2n-2j_2+1) +2(j_3-j_2)(2n-2j_3+1),\\
    \dim_{\mathrm{CR}}M&=2n+2j_3-7,\\
    \mathrm{co}\dim_{\mathrm{CR}}M&=4\big((n-j_3)(j_3-1)
    +(j_3-j_2)(j_2-1)+1\big).
  \end{align*}
  The space of Levi
  forms of $M$ has dimension \[ (2j_2-3)(2j_3-2j_2+1)
  +(2j_3-2j_2-1)(2n-2j_3+2j_2-1) +(2n-2j_3)(2j_3-2j_2),\] and all of
  them, in a suitable basis (as described in \cite{AMN}), have all
  diagonal entries equal to zero. Hence $M$ is trace pseudoconcave,
  and has not kind $2$. It can be checked that $M$ has kind $3$, thus
  $M$ satisfies the hypotheses of Theorem~\ref{main} and germs of CR
  distributions extend holomorphically to a full neighborhood.

  With the notation of \cite{AMN} $M$ is the minimal orbit of the
  simple Lie group of type $\mathrm{A\,II}_{2n-1}$ associated to the
  parabolic subalgebra
  $\mathfrak{q}=\mathfrak{q}_{\{\alpha_{1},\alpha_{2j_2-1},\alpha_{2j_3-1}\}}$. By
  duality, a completely analogous statement holds for the case $1\leq
  j_1<j_2<j_3=n$.
\end{example}

\begin{example}\rm
  Similar to the previous one is the case of pairs of
  complex flag manifolds, homogeneous for the action of
  $\mathrm{SL}(2n,\mathbb{C})$, and compact generic CR submanifolds,
  homogeneous for the action of $\mathrm{SL}(n,\mathbb{H})$:
  \begin{equation*}
  \begin{cases}
      Y=\big\{\ell_{2j_1-1}\subset\ell_{2j_2-1}
        \subset\mathbb{C}^{2n}\mid
        \dim\ell_{2j_k-1}=2j_k-1,\; k=1,2\big\}, \\[4pt]
      M=\left\{(\ell_{2j_1-1},\ell_{2j_2-1})\in Y\left|
        \begin{aligned}&\dim(\ell_{2j_k-1}\cap\mathbb H\ell_{2j_k-1})=2j_k-2,
        &&k=1,2\\
        &\;\mathbb H\ell_{2j_1-1}\subset\ell_{2j_{2}-1}\end{aligned}
        \right.\right\},\\[4pt]
      \text{for}\quad 1< j_1<j_2<n.
    \end{cases}
  \end{equation*}
  Also in this case, $M$ has kind $3$ and is trace pseudoconcave.
\end{example}

\begin{example}\rm
  Let $G^{\mathbb{C}}$ be the connected and simply connected simple
  complex Lie group of type $\mathrm F_4$, and let
  $\mathfrak{g}^\mathrm{C}$ its Lie algebra. Fix a Cartan subalgebra
  $\mathfrak{h}^{\mathbb{C}}$ of $\mathfrak{g}^\mathrm{C}$ and a
  system of simple roots $\mathcal B=\{\alpha_j\}_{1\leq j\leq4}$ of
  the root system $\mathcal R=\mathcal
  R(\mathfrak{g}^\mathrm{C},\mathfrak{h}^{\mathbb{C}})$ (we use the
  root numbering scheme of Bourbaki, see \cite{AMN}, Appendix). Let
  $\{\omega_j\}_{1\leq j\leq4}$ be the set of fundamental weights dual
  to $\mathcal{B}$ and let $V=V_{\omega_2}$ be the fundamental
  representation of $G^{\mathbb{C}}$ with highest weight $\omega_2$.
  Let $Y$ be the
  $G^{\mathbb{C}}$-orbit, in $\mathbb{P}(V)$, of
  the highest weight root space.
  With the notation of \cite{AMN}, $Y$ is the complex flag
  manifold of $G^{\mathbb{C}}$ consisting of parabolic subalgebras of
  $\mathfrak{g}^{\mathbb{C}}$ that are
  conjugate to $\mathfrak{q}_{\{\alpha_2\}}$.

  Inside $Y$ consider the minimal orbit $M$ of the real form of
  $G^\mathbb{C}$ of type $\mathrm{F\,II}$ (also denoted by
  $F_{4(-20)}$).  It is a CR manifold of CR dimension $9$ and
  CR codimension $11$.  By \cite{AMN}, Theorem~9.1, it is of finite
  type, and by \cite{AMN}, Theorem~13.5, it is trace
  pseudoconcave. Direct computation shows that $G_2$ has codimension
  $4$ and $G_3$ is the whole tangent space of $M$, thus $M$ has kind
  $3$. Hence also in this case we have holomorphic extension for germs
  of CR distributions.
\end{example}


\begin{thebibliography}{CLN1985}


\bibitem[AMN]{AMN}
{\sc Altomani, A., Medori, C., Nacinovich, M.}:
{\sl The CR structure of minimal orbits in complex flag manifolds}.
J.~Lie Theory {\bf 16} (2006), 483--530.


\bibitem[BR]{BR}
{\sc Baouendi, M.~S.~, Rothschild, L.~P.}:
{\sl Cauchy Riemann functions on manifolds of higher codimension
in complex space}.
Invent.~Math.~{\bf 101} (1990), 45-56.

\bibitem[BT1]{BT1}
{\sc Baouendi, M.~S.~, Treves, F.}:
{\sl A property of the functions and distributions annihilated by a
locally integrable system of complex vector fields}.
Ann.~Math.~{\bf 113} (1981), 387-421.

\bibitem[BT2]{BT2}
{\sc Baouendi, M.~S.~, Treves, F.}:
{\sl About the holomorphic extension of CR functions on
real hypersurfaces in complex space}.
Duke Math. J. {\bf 51} (1984),  no. 1, 77--107.

\bibitem[B]{B}
{\sc Bishop, E.}:
{\sl Differentiable manifolds in complex Euclidean space}.
Duke Math.~J.~{\bf 32} (1965), 1-21.

\bibitem[Bo]{Bo}
{\sc Boggess A.}:
{\sl CR extendability near a point where the first Levi form vanishes}.
Duke Math.~J.~{\bf 48} (1981), 665-684.

\bibitem[BP]{BP}
{\sc Boggess A., Polking J.}:
{\sl Holomorphic extensions of $CR$ functions}.
Duke Math. J. {\bf 49} (1982), 757-784.


\bibitem[DS]{DS}
{\sc Dinh, T.~C., Sarkis, F.~}:
{\sl Wedge removability of metrically thin sets and application to the CR
meromorphic extension}.
Math.~Z.~{\bf 238} (2001), 639-653.

\bibitem[DH]{DH}
{\sc Dolbeault, P., Henkin, G.}:
{\sl Cha\^{\i}nes holomorphes de bord donn\'e dans un ouvert q-concave de $\mathbb{C}P^n$}.
Bull.~Soc.~Math.~France {\bf 125} (1997), 383-445.

\bibitem[FR]{FR}
{\sc Fornaess, J.~E., Rea, C.}:
{\sl Local holomorphic extendability and nonextendability
of CR-functions on smooth boundaries}.
Ann. Scuola Norm. Sup. Pisa Cl. Sci. (4) {\bf 12} (1985), no. 3, 491--502.

\bibitem[HL]{HL}
{\sc Harvey, R., Lawson, B.}:
{\sl On boundaries of complex analytic varieties, Part II}.
Ann.~Math. {\bf 106} (1977), 213-238.

\bibitem[HN1]{HN1}
{\sc Hill, C.~D., Nacinovich, M.}:
{\sl A weak pseudoconcavity condition for abstract almost CR manifolds}.
Invent.~Math.~{\bf 142} (2000), 251-283.

\bibitem[HN2]{HN2}
{\sc Hill, C.~D., Nacinovich, M.}:
{\sl Fields of CR meromorphic functions}.
Rend. Sem. Mat. Univ. Padova {\bf 111} (2004), 179-204.

\bibitem[HN3]{HN3}
{\sc Hill, C. D., Nacinovich, M.}:
{\sl Elementary pseudoconcavity and fields of CR meromorphic functions}.
Rend. Sem. Mat. Univ. Padova 113 (2005), 99-115.

\bibitem[HN4]{HN4}
{\sc Hill, C. D., Nacinovich, M.}:
{\sl Conormal suspensions of differential complexes}.
J. Geom. Anal. 10, 3 (2000), 496-537.


\bibitem[HT]{HT}
{\sc Hill, C.~D., Taiani, G.}:
{\sl Families of analytic discs in $\C^n$ with boundaries on a
prescribed CR submanifold}.
Ann.~Sc.~Norm.~Pisa {\bf 5} (1978), 327-380.


\bibitem[J]{J}
{\sc J\"oricke, B.}:
{\sl Deformation of CR-manifolds, minimal points and CR-manifolds with the microlocal
analytic extension property}.
J.~Geom.~Anal.~{\bf 6} (1996), 555-611.

\bibitem[MN1]{MN1}
{\sc Medori, C., Nacinovich, M.}:
{\sl Levi-Tanaka algebras and homogeneous CR manifolds}.
Compositio Mathematica {\bf 109} (1997), 195–-250.

\bibitem[MN2]{MN2}
{\sc Medori, C., Nacinovich, M.}:
{\sl Classification of semisimple Levi-Tanaka algebras}.
Ann.~Mat.~Pura Appl.~{\bf CLXXIV} (1998), 285–-349.

\bibitem[MN3]{MN3}
{\sc Medori, C., Nacinovich, M.}:
{\sl Complete nondegenerate locally standard CR manifolds}.
Math.~Ann.~{\bf 317} (2000), 509-526.

\bibitem[MN4]{MN4}
{\sc Medori, C., Nacinovich, M.}:
{\sl Algebras of infinitesimal CR automorphisms}.
J.~Algebra {\bf 287} (2005), 234-274.

\bibitem[M]{M}
{\sc Merker, J.}:
{\sl Global minimality of generic manifolds and holomorphic extendibility
of CR functions}.
IMRN {\bf 8} (1994), 329-342.

\bibitem[MP1]{MP1}
{\sc Merker, J., Porten, E.}:
{\sl On the local meromorphic extension of CR meromorphic mappings}.
Ann.~Polon.~Math.~{\bf 70} (1998), 163-193.

\bibitem[MP2]{MP2}
{\sc Merker, J., Porten, E.}:
{\sl Metrically thin singularities of integrable CR functions}.
Internat.~J.~Math.~{\bf 11} (2000), 857-872.

\bibitem[MP3]{MP3}
{\sc Merker, J., Porten, E.}:
{\sl Characteristic foliations on maximally real submanifolds of $\C^n$
and removable singularities for CR functions}.
IMRP Volume {\bf 2006} (2006), Article ID 72069,
p.1-131 .

\bibitem[MP4]{MP4}
{\sc Merker, J., Porten, E.}:
{\sl Holomorphic Extension of CR Functions, Envelopes of Holomorphy, and Removable
Singularities}.
IMRS Volume {\bf 2006} (2006),
Article ID 28925, p.1-286.

\bibitem[NV]{NV}
{\sc Nacinovich, M., Valli, G.}:
{\sl Tangential Cauchy-Riemann complexes on distributions}.
Ann. Mat. Pura Appl. {\bf 146} (1987), 123-160.

\bibitem[S]{S}
{\sc Sj\"ostrand, J.}:
{\sl The FBI transform for CR submanifolds of $\C^n$}.
Pr\'epublications Math\'ematiques Orsay 1982.

\bibitem[T1]{T1}
{\sc Tr\'epreau, J.-M.}:
{\sl Sur le prolongement holomorphe des fonctions CR d\'efinies
sur une hypersurface reelle de classe $C^2$}.
Invent.~Math.~{\bf 83} (1986), 583-592.

\bibitem[T2]{T2}
{\sc Tr\'epreau, J.-M.}:
{\sl Sur la propagation des singularit\'es dans les vari\'et\'es CR}.
Bull.~Soc. Math.~Fr.~{\bf 118} (1990), 403-450.

\bibitem[T]{T}
\newblock{F.~Treves.}
\newblock{\it Hypoanalytic Structures: Local Theory.}
\newblock{Princeton Univ. Press 1992.}

\bibitem[Tu1]{Tu1}
{\sc Tumanov, A.~E.}:
{\sl Extension of CR-functions into a wedge from a manifold of finite type}.
Math.~USSR Sb.~{\bf 64} (1989), 129-140.

\bibitem[Tu2]{Tu2}
{\sc Tumanov, A.~E.}:
{\sl On the propagation of extendability of CR functions}.
Complex Analysis and Geometry, Proc.~ of the Conference in Trento 1995,
Lecture Notes in Pure and Appl.~Math.~{\bf 173}, Dekker 1996, pp.479-498.

\end{thebibliography}
\end{document}